\documentclass[12pt,reqno]{amsart}

\usepackage[dvips]{graphicx}
\usepackage[arrow,matrix,curve]{xy}

\usepackage{amssymb, latexsym, amsmath, amscd, array, 
}

\newtheorem{theorem}{Theorem}[section]
\newtheorem{lemma}[theorem]{Lemma}
\newtheorem{proposition}[theorem]{Proposition}
\newtheorem{corollary}[theorem]{Corollary}
\newtheorem{question}[theorem]{Question}

\theoremstyle{definition}
\newtheorem{definition}[theorem]{Definition}
\newtheorem{example}[theorem]{Example}
\newtheorem{remark}[theorem]{Remark}

\numberwithin{equation}{section}
\numberwithin{figure}{section}
\numberwithin{table}{section}

\DeclareMathOperator{\rhp2}{{\mathcal P}}

\DeclareMathOperator{\PD}{{\rm PD}}

\DeclareMathOperator{\sysh}{{\rm sysh}}

\DeclareMathOperator{\stsys}{{\rm stsys}}
\DeclareMathOperator{\Confsys}{{\rm Confsys}}
\DeclareMathOperator{\confsys}{{\rm confsys}}

\DeclareMathOperator{\trace}{{\rm trace}}

\DeclareMathOperator{\area}{{\rm area}}
\DeclareMathOperator{\vol}{{\rm vol}}

\newcommand\Spin{{\rm Spin}}
\newcommand\SO{{\operatorname{SO}}}
\newcommand\form {{f}}

\newcommand\kraines{{\kappa_{\HP}^{\phantom{I}}}}
\newcommand\krainesl{{\kappa}}

\newcommand\ie {{\it i.e.\ }}

\newcommand\cf {\hbox{\it cf.\ }}

\newcommand\gmetric{{\mathcal G}}

\newcommand{\syspi} {{{\rm sys}\pi}}

\newcommand{\surface}{\Sigma}

\newcommand\SR{{\rm SR}}

\newcommand\cay {\omega_{\rm Ca}^{\phantom{I}}}
\newcommand\cayp {\omega_{\rm Ca}^{||}}

\newcommand \cqfd{\unskip\kern 6pt\penalty 500
\raise -2pt\hbox{\vrule\vbox to10pt{\hrule width 4pt
\vfill\hrule}\vrule}\par}

\def\adots{\mathinner{\mkern2mu\raise1pt\hbox{.}
\mkern3mu\raise4pt\hbox{.}\mkern1mu\raise7pt\hbox{.}}}
\def\hfl#1{\frac{\buildrel{#1}}{{\hbox to 12mm{\rightarrowfill}}}}
  
\def \\R^n \times \R^n
\rightarrow \R{\mathop{\R^n \times \R^n
\rightarrow \R}}

\newcommand\C{{\mathbb {C}}}

\newcommand\HH {{\mathbb H}}
\newcommand\HP {{{\mathbb H}{\mathbb P}}}

\newcommand\N {{\mathbb N}}

\newcommand\PP {{\mathbb {P}}}

\newcommand\Q{{\mathbb {Q}}} \newcommand\R {{\mathbb R}}
\newcommand\RR {{\mathbb R}} \newcommand\RP {{\mathbb R}{\mathbb P}}
 \newcommand\Z {{\mathbb Z}} \newcommand\W
{{\rm Wirt}} \newcommand\Wirt {{\rm Wirt}}

\newcommand\CP {{\mathbb C\mathbb P}}

\newcommand{\cptwo}{\mathbb C\mathbb P^2}
\newcommand{\hypo}{{\mathcal J}}

\newcommand\fake {{homology}}

\long\def\forget#1\forgotten{} %

\long\def\forgett#1\forgottent{} %

\def\circ{\mathchoice%
 {\mathrel{\raise 1pt\hbox{$\scriptstyle\mathchar"020E$}}}
 {\mathrel{\raise 1pt\hbox{$\scriptstyle\mathchar"020E$}}}
 {\mathrel{\raise 1pt\hbox{$\scriptscriptstyle\mathchar"020E$}}}
 {}
}

\newcommand{\nc}{\newcommand} \nc{\on}{\operatorname}

\nc{\df}{\on{\it df}}

\nc{\conf}{\on{conf}}

\nc{\spt}{\on{spt}}
\nc{\norm}[1]{\| #1 \|}
\nc{\parallelleer}{\norm{\ }} 
\nc{\parallelh}{\norm h} 
\nc{\parallelk}{\norm k} 
\nc{\parallelx}{\norm x} 
\nc{\parallelhrr}{\norm {h_\RR}} 
\nc{\parallelom}{\norm \omega} 
\nc{\parallelomij}{\norm {\omega_{i_j}}} 
\nc{\parallelomx}{\norm {\omega_{x}}} 
\nc{\parallelpi}{\norm \pi} 
\nc{\parallelalf}{\norm \alpha} 
\nc{\parallelalfs}{\norm {\alpha_s}} 
\nc{\parallelalfi}{\norm {\alpha_i}} 
\nc{\parallelalfij}{\norm {\alpha_{i_j}}} 
\nc{\parallelbeta}{\norm \beta} 
\nc{\parallelbetat}{\norm {\beta_t}} 
\nc{\parallelhcapalf}{\norm {h \cap \alpha}} 
\nc{\parallelPDralf}{\norm {\PD_\RR(\alpha)}} 
\nc{\strichleer}{| \  |}
\nc{\NN}{\mathbb N}

\nc{\rr}{\mbox{$\scriptstyle\mathbb R$}}

\nc{\dF}{{\it dF}}

\nc{\DF}{{\it DF}}

\nc{\ds}{{\it ds}}

\nc{\dvol}{{\it dvol}}

\nc{\grad}{{\rm grad}}
\nc{\strichw}{\|\omega\|}
\nc{\strichwx}{|\omega_x|}
\nc{\Hess}{{\rm Hess}}

\begin{document}

\title{$E_7$, Wirtinger inequalities, Cayley~$4$-form, and homotopy}

\author[V.~Bangert]{Victor Bangert$^*$}

\address{ Mathematisches Institut, Universit\"at Freiburg,
Eckerstr.~1, 79104 Freiburg, Germany}

\email{bangert@mathematik.uni-freiburg.de}

\thanks{$^{*}$Partially Supported by DFG-Forschergruppe `Nonlinear
Partial Differential Equations: Theoretical and Numerical Analysis'}

\author[M.~Katz]{Mikhail G. Katz$^{**}$}

\address{Department of Mathematics, Bar Ilan University, Ramat Gan
52900 Israel}

\email{katzmik@macs.biu.ac.il}

\thanks{$^{**}$Supported by the Israel Science Foundation (grants
84/03 and 1294/06) and the BSF (grant 2006393)}

\author[S.~Shnider]{Steven Shnider}

\address{Department of Mathematics, Bar Ilan University, Ramat Gan
52900 Israel}

\email{shnider@macs.biu.ac.il}

\author[S.~Weinberger]{Shmuel Weinberger$^{***}$}

\address{Department of Mathematics, University of Chicago, Chicago, IL
60637}

\email{shmuel@math.uchicago.edu}

\thanks{$^{***}$Partially supported by NSF grant DMS 0504721 and the
BSF (grant 2006393)}



\date{\today}

\begin{abstract}
We study optimal curvature-free inequalities of the type discovered by
C.~Loewner and M.~Gromov, using a generalisation of the Wirtinger
inequality for the comass.  Using a model for the classifying
space~$BS^3$ built inductively out of~$BS^1$, we prove that the
symmetric metrics of certain two-point homogeneous manifolds turn out
not to be the systolically optimal metrics on those manifolds.  We
point out the unexpected role played by the exceptional Lie
algebra~$E_7$ in systolic geometry, via the calculation of Wirtinger
constants.  Using a technique of pullback with controlled systolic
ratio, we calculate the optimal systolic ratio of the quaternionic
projective plane, modulo the existence of a Joyce manifold with
Spin(7) holonomy and unit middle-dimensional Betti number.
\end{abstract}

\maketitle
\tableofcontents

\section{Inequalities of Pu and Gromov}

The present text deals with systolic inequalities for the projective
spaces over the division algebras~$\R$,~$\C$, and~$\HH$.

In 1952, P.M. Pu \cite{Pu} proved that the least length,
denoted~$\syspi_1$, of a noncontractible loop of a Riemannian
metric~$\gmetric$ on the real projective plane~$\RP^2$, satisfies the
optimal inequality
\[
\syspi_1(\RP^2, \gmetric)^2 \leq \tfrac{\pi}{2} \area(\RP^2,\gmetric).
\]
Pu's bound is attained by a round metric, \ie one of constant Gaussian
curvature.  This inequality extends the ideas of C.~Loewner, who
proved an analogous inequality for the torus in a graduate course at
Syracuse University in 1949, thereby obtaining the first result in
systolic geometry, \cf \cite{SGT}.

Defining the optimal systolic ratio~$\SR(\surface)$ of a
surface~$\surface$ as the supremum
\begin{equation}
\label{11b}
\SR(\surface) = \sup_\gmetric \left\{
\left. \frac{\syspi_1(\gmetric)^2} {\area (\gmetric)} \right|
~\gmetric~{\rm Riemannian~metric~on}~\surface \right\},
\end{equation}
we can restate Pu's inequality as the calculation of the value
\[
\SR(\RP^2)= \tfrac{\pi}{2},
\]
the supremum being attained by a round metric.

One similarly defines a homology systole, denoted~$\sysh_1$, by
minimizing over loops in~$\surface$ which are not nullhomologous.  One
has~$\syspi_1(\surface) \leq \sysh_1 (\surface)$.  For orientable
surfaces, one has the identity
\begin{equation}
\label{12b}
\sysh_1(\surface)= \lambda_1 \left( H_1^{\phantom{I}}(\surface,\Z),
\|\;\| \right),
\end{equation}
where~$\|\;\|$ is the stable norm in homology (see Section~\ref{gro}),
while~$\lambda_1$ is the first successive minimum of the normed
lattice.  In other words, the homology systole and the
stable~$1$-systole (see below) coincide in this case (and more
generally in codimension~$1$).  Thus, the homology 1-systole is the
least stable norm of an integral~$1$-homology class of infinite order.

Therefore, either the homology~$k$-systole or the stable~$k$-systole
can be thought of as a higher-dimensional generalisation of the
1-systole of surfaces.  It has been known for over a decade that the
homology systoles do not satisfy systolic inequalities; see \cite{Ka}
where the case of the products of spheres~$S^k\times S^k$ was treated.
Homology systoles will not be used in the present text.

For a higher dimensional manifold~$M^{2k}$, the appropriate
middle-dimensional invariant is therefore the stable~$k$-systole
$\stsys_k$, defined as follows.  Let~$H_k(M, \Z)_\R^{\phantom{I}}$ be
the image of the integral lattice in real~$k$-dimensional homology
of~$M$.  The~$k$-Jacobi torus~$J_k M$ is the quotient
\begin{equation}
\label{13}
J_k M= H_k(M,\R)/H_k(M,\Z)_R.
\end{equation}
We set
\begin{equation}
\label{11z}
\stsys_k(\gmetric)=\lambda_1\left( H_k(M, \Z)_\R^{\phantom{I}}, \|\;\|
\right),
\end{equation}
where~$\|\;\|$ is the stable norm in homology, while~$\lambda_1$ is
the first successive minimum of the normed lattice.  In other words,
the stable~$k$-systole is the least stable norm of an
integral~$k$-homology class of infinite order.  A detailed definition
of the stable norm appears in Section~\ref{gro}.

By analogy with \eqref{11b}, one defines the optimal
middle-dimensional stable systolic ratio,~$\SR_k (M^{2k})$, by setting
\[
\SR_k(M) = \sup_\gmetric \frac{\stsys_k(\gmetric)^2} {\vol_{2k}
(\gmetric)} \,,
\]
where the supremum is over all Riemannian metrics~$\gmetric$ on~$M$.

In 1981, M. Gromov \cite{Gr0} proved an inequality analogous to Pu's,
for the complex projective plane~$\CP^2$.  Namely, he evaluated the
optimal stable systolic ratio of~$\CP^2$, which turns out to be
\[
\SR_2(\CP^2) = 2,
\]
where, similarly to the real case, the implied optimal bound is
attained by the symmetric metric, \ie the Fubini-Study metric.  In
fact, Gromov proved a more general optimal inequality.

\begin{theorem}[M. Gromov]
Every metric~$\gmetric$ on the complex projective space satisfies the
inequality
\begin{equation}
\label{12}
\stsys_2(\CP^n, \gmetric)^n \leq n! \vol_{2n}(\CP^n,\gmetric).
\end{equation}
\end{theorem}

Here~$\stsys_2$ is still defined by formula~\eqref{11z} with~$k=2$,
and we set~$M=\CP^n$.

A quaternionic analogue of the inequalities of Pu and Gromov was
widely expected to hold.  Namely, the symmetric metric on the
quaternionic projective plane~$\HP^2$ gives a ratio equal
to~$\tfrac{10}{3}$, calculated by a calibration argument in
Section~\ref{1215}, following the approach of \cite{Be1}.  It was
widely believed that the {\em optimal\/} systolic ratio~$\SR_4(\HP^2)$
equals~$\tfrac{10}{3}$, as well.  See also \cite[Section~4]{Gr2} and
\cite[Remark~4.37, p.~262]{Gr3} or~\cite{Gr4}.  Contrary to
expectation, we prove the following theorem.

\begin{theorem}
\label{11}
The quaternionic projective space~$\HP^{2n}$ and the complex
projective space~$\CP^{4n}$ have a common optimal middle dimensional
stable systolic ratio:~$\SR_{4n}(\HP^{2n}) = \SR_{4n}(\CP^{4n})$.
\end{theorem}

Theorem~\ref{11} is proved in Section~\ref{sevenb}.  The Fubini-Study
metric gives a middle-dimensional ratio equal to~$(4n)!/((2n)!)^2$ for
the complex projective~$4n$-space.  For instance, the symmetric metric
of~$\CP^4$ gives a ratio of~$6$.  The symmetric metric on~$\HP^{2n}$
has a systolic ratio of~$(4n+1)!/((2n+1)!)^2$, \cf \cite{Be1}.  Since
\[
(4n+1)!/((2n+1)!)^2 < (4n)!/((2n)!)^2,
\]
we obtain the following corollary.

\begin{corollary}
The symmetric metric on~$\HP^{2n}$ is not systolically optimal.
\end{corollary}

We also estimate the common value of the optimal systolic ratio in the
first interesting case, as follows.

\begin{proposition}
\label{11y}
The common value of the optimal ratio for~$\HP^2$ and~$\CP^4$ lies in
the following interval:
\begin{equation}
\label{15}
6 \leq \SR_4(\HP^2) = \SR_4(\CP^4) \leq 14.
\end{equation}
\end{proposition}

The constant~$14$ which appears above as the upper bound for the
optimal ratio, is twice the dimension of the Cartan subalgebra of the
exceptional Lie algebra~$E_7$, reflected in our title.  More
specifically, the relevant ingredient is that every self-dual~$4$-form
admits a decomposition into at most~$14$ decomposable (simple) terms
with respect to a suitable orthonormal basis, \cf proof of
Proposition~\ref{45}.

Note that quaternion algebras and congruence subgroups of arithmetic
groups were used in \cite{KSV} to study asymptotic behavior of the
systole of Riemann surfaces.  It was pointed out by a referee that for
the first time in the history of systolic geometry, Lie algebra theory
has been used in the field.

We don't know of any techniques for constructing metrics on~$\CP^4$
with ratio greater than the value~$6$, attained by the Fubini-Study
metric.  Meanwhile, an analogue of Gromov's proof for~$\CP^2$ only
gives an upper bound of~$14$.  This is due to the fact that the
Cayley~$4$-form~$\cay$, \cf \cite{Be1,HL}, has a higher Wirtinger
constant than does the Kahler~$4$-form (\ie the square of the standard
symplectic~$2$-form).  Nonetheless, we expect that the resulting
inequality is optimal, \ie that the common value of the optimal
systolic ratio of~$\HP^2$ and~$\CP^4$ is, in fact, equal to~$14$.  The
evidence for this is the following theorem, which should give an idea
of the level of difficulty involved in evaluating the optimal ratio in
the quaternionic case, as compared to Pu's and Gromov's calculations.
Joyce manifolds \cite{Jo0} are discussed in Section~\ref{ten}.

\begin{theorem}
\label{14a}
If there exists a compact Joyce manifold~$\hypo$ with~$\Spin(7)$
holonomy and with~$b_4(\hypo)=1$, then the common value of the middle
dimensional optimal systolic ratio of~$\HP^2$ and~$\CP^4$ equals~$14$.
\end{theorem}

A smooth Joyce manifold with middle Betti number~$1$ would necessarily
be rigid.  Thus it cannot be obtained by any known techniques, relying
as they do on deforming the manifold until it decays into something
simpler.  On the other hand, by relaxing the hypothesis of smoothness
to, say, that of a PD(4) space, such a mildly singular Joyce space may
be obtainable as a suitable quotient of an~$8$-torus, and may be
sufficient for the purposes of calculating the systolic ratio in this
dimension.

\begin{corollary}
If there exists a compact Joyce manifold~$\hypo$ with~$\Spin(7)$
holonomy and with~$b_4(\hypo)=1$, then the symmetric metric on~$\CP^4$
is not systolically optimal.
\end{corollary}

If one were to give a synopsis of the history of the application of
homotopy techiques in systolic geometry, one would have to start with
D.~Epstein's work \cite{Ep} on the degree of a map in the 1960's,
continue with A.~Wright's work \cite{Wr} on monotone mappings in the
1970's, then go on to developments in real semi-algebraic geometry
which indicated that an arbitrary map can be homotoped to have good
algebraic structure by M. Coste and others \cite{BCR}, in the 1980's.

M. Gromov, in his 1983 paper \cite{Gr1}, goes out of the category of
manifolds in order to prove the main isoperimetric inequality relating
the volume of a manifold, to its filling volume.  Namely, the cutting
and pasting constructions in the proof of the main isoperimetric
inequality involve objects more general than manifolds.

In the 1992 paper in Izvestia by I. Babenko \cite{Bab1}, his Lemma 8.4
is perhaps the place where a specific homotopy theoretic technique was
first applied to systoles.  Namely, this technique derives
systolically interesting consequences from the existence of maps from
manifolds to simplicial complexes, by pullback of metrics.  This work
shows how the triangulation of a map~$f$, based upon the earlier
results mentioned above, can help answer systolic questions, such as
proving a converse to Gromov's central result of 1983.  What is
involved, roughly, is the possibility of pulling back metrics by~$f$,
once the map has been deformed to be sufficiently nice (in particular,
real semialgebraic).

In 1992-1993, Gromov realized that a suitable oblique~$\Z$ action on
the product~$S^3 \times \R$ gives a counterexample to
a~$(1,3)$-systolic inequality on the product~$S^1 \times S^3$.  This
example was described by M.~Berger~\cite{Be5}, who sketched also
Gromov's ideas toward constructing further examples of systolic
freedom.

In 1995, metric simplicial complexes were used \cite{Ka} to prove the
systolic freedom of the manifold~$S^n \times S^n$.  In this paper, a
polyhedron~$P$ is defined in equation (3.1).  It is exploited in an
essential way in an argument in the last paragraph on page 202, in the
proof of Proposition~3.3.

Thus, we will exploit a map of classifying spaces~$BS^1 \to BS^3$ so
as to relate the systolic ratios of the quaternionic projective space
and the complex projective space.  We similarly relate the
quaternionic projective space and a hypothetical Joyce manifold
(with~$\Spin_7$ holonomy) with~$b_4=1$, relying upon a result by
H. Shiga in rational homotopy theory.

An interesting related axiomatisation (in the case of 1-systoles) is
proposed by M. Brunnbauer \cite{Bru}, who proves that the optimal
systolic constant only depends on the image of the fundamental class
in the classifying space of the fundamental group, generalizing
earlier results of I. Babenko.  For background systolic material, see
\cite{Gr1, Ka, BK2, KL, SGT}.

In Section~\ref{gro}, we present Gromov's proof of the optimal stable
2-systolic inequality~\eqref{12} for the complex projective
space~$\CP^n$, \cf \cite[Theorem 4.36]{Gr3}, based on the cup product
decomposition of its fundamental class.  The proof relies upon the
Wirtinger inequality, proved in Section~\ref{fed} following H.~Federer
\cite{Fe1}.  In Section~\ref{1215}, we analyze the symmetric metric on
the quaternionic projective plane from the systolic viewpoint.  A
general framework for Wirtinger-type inequalities is proposed in
Section~\ref{e7}.

A homotopy equivalence between~$\HP^n$ and a suitable CW complex built
out of~$\CP^{2n}$ is constructed in Section~\ref{beegee} using a map
$BS^1 \to B S^3$.  Section~\ref{sevenb} exploits such a homotopy
equivalence to build systolically interesting metrics.
Section~\ref{eight} contains some explicit formulas in the context of
the Kraines form and the Cayley form~$\cay$.  Section~\ref{nine}
presents a Lie-theoretic analysis of~$4$-forms on~$\R^8$, using an
idea of G.~Hunt.  Theorem~\ref{14a} is proved in Section~\ref{ten}.
Related results on the Hopf invariant and Whitehead products are
discussed in Section~\ref{eleven}.

\section{Federer's proof of Wirtinger inequality}
\label{fed}

Following H.~Federer \cite[p.~40]{Fe1}, we prove an optimal upper
bound for the comass norm~$\| \; \|$, \cf Definition~\ref{34}, of the
exterior powers of a~$2$-form.

Recall that an exterior form is called {\em simple\/} (or {\em
decomposable}) if it can be expressed as a wedge product of~$1$-forms.
The comass norm for a simple~$k$-form coincides with the natural
Euclidean norm on~$k$-forms.  In general, the comass is defined as
follows.

\begin{definition}
\label{34}
The comass of an exterior~$k$-form is its maximal value on a~$k$-tuple
of unit vectors.
\end{definition}

Let~$V$ be a vector space over~$\C$.  Let~$H=H(v,w)$ be a Hermitian
product on~$V$, with real part~$v\cdot w$, and imaginary
part~$A=A(v,w)$, where~$A\in \bigwedge^2 V$, the second exterior power
of~$V$.  Here we adopt the convention that~$H$ is complex linear in
the {\em second\/} variable.

\begin{example}
\label{1211z}
Let~$Z_1,\ldots, Z_\nu \in \bigwedge^1(\C^\nu, \C)$ be the coordinate
functions~ in~$\C^\nu$.  We then have the standard
(symplectic)~$2$-form, denoted~$A \in \bigwedge^2 (\C^\nu, \C)$, given
by
\[
A= \tfrac{i}{2} \sum_{j=1}^\nu Z_j\wedge \bar Z_j .
\]
\end{example}

\begin{lemma}
\label{1213}
The comass of the standard symplectic form~$A$ satisfies~$\| A \| =1$.
\end{lemma}

\begin{proof}
We can set~$\xi= v \wedge w$, where~$v$ and~$w$ are orthonormal.  We
have~$H(v,w)= iA(v,w)$, hence
\begin{equation}
\label{21}
\langle \xi, A \rangle = A(v,w) = H(iv, w) = (iv) \cdot w \leq 1
\end{equation}
by the Cauchy-Schwarz inequality; equality holds if and only if one
has~$~iv=w$.
\end{proof}

\begin{remark}
R.~Harvey and H.~B.~Lawson \cite{HL} provide a similar argument for
the Cayley~$4$-form~$\cay$.  They realize~$\cay$ as the real part of a
suitable multiple vector product on~$\R^8$, defined in terms of the
(non-associative) octonion multiplication, to calculate the comass
of~$\cay$, \cf Proposition~\ref{83b}.
\end{remark}

\begin{proposition}[Wirtinger inequality]
\label{23}
Let~$\mu \geq 1$.  If~$\xi\in \bigwedge_{2\mu} V$ and~$\xi$ is simple,
then
\[
\langle \xi, A^\mu \rangle \leq \mu!\; |\xi|;
\]
equality holds if and only if there exist elements~$v_1,\ldots,v_\mu
\in V$ such that
\[
\xi= v_1 \wedge (iv_1) \wedge \cdots \wedge v_\mu \wedge (iv_\mu).
\]
Consequently,~$\| A^\mu \| = \mu !$
\end{proposition}

\begin{proof}
The main idea is that in real dimension~$2\mu$, every~$2$-form is
either simple, or splits into a sum of at most~$\mu$ orthogonal simple
pieces.

We assume that~$|\xi|=1$, where~$|\;|$ is the natural Euclidean norm
in~$\bigwedge_{2\mu}V$.  The case~$\mu=1$ was treated in
Lemma~\ref{1213}.

In the general case~$\mu \geq 1$, we consider the~$2\mu$ dimensional
subspace~$T$ associated with~$\xi$.  Let~$f: T \to V$ be the inclusion
map, and consider the pullback 2-form~$(\wedge^2 f) A \in \bigwedge^2
T$.  Next, we orthogonally diagonalize the skew-symmetric 2-form, \ie
decompose it into~$2\times 2$ diagonal blocks.  Thus, we can choose
dual orthonormal bases~$e_1,\ldots, e_{2\mu}$ of~$T$
and~$\omega_1,\ldots, \omega_{2\mu}$ of~$\bigwedge^1 T$, and
nonnegative numbers~$\lambda_1, \ldots, \lambda_\mu$, so that
\begin{equation}
\label{1211a}
(\wedge^2 f) A = \sum_{j=1}^\mu \lambda_j \left( \omega_{2j-1} \wedge
\omega_{2j} \right) .
\end{equation}
By Lemma~\ref{1213}, we have
\begin{equation}
\label{22}
\lambda_j = A(e_{2j-1}, e_{2j}) \leq \|A\| = 1
\end{equation}
for each~$j$.  Noting that~$\xi= \epsilon e_1 \wedge \cdots \wedge
e_{2\mu}$ with~$\epsilon= \pm 1$, we compute
\[
\left( \wedge^{2\mu} f \right) A^\mu = \mu! \lambda_1 \ldots
\lambda_\mu \omega_1 \wedge \cdots \wedge \omega_{2\mu},
\]
and therefore
\begin{equation}
\label{1212z}
\langle \xi, A^\mu \rangle = \epsilon \mu!\; \lambda_1 \ldots
\lambda_\mu \leq \mu!
\end{equation}
Note that equality occurs in \eqref{1212z} if and only if~$\epsilon
=1$ and~$\lambda_j=1$.  Applying the proof of Lemma~\ref{1213}, we
conclude that~$e_{2j} = i e_{2j-1}$, for each~$j$.
\end{proof}

\begin{corollary}
\label{1214}
Every real~$2$-form~$A$ satisfies the comass bound
\begin{equation}
\label{1212}
\| A^\mu \| \leq \mu! \| A\|^\mu.
\end{equation}
\end{corollary}

\begin{proof}
An inspection of the proof Proposition~\ref{23} reveals that the
orthogonal diagonalisation argument, \cf \eqref{22}, applies to an
arbitrary~$2$-form~$A$ with comass~$\| A \| = 1$.
\end{proof}

\begin{lemma}
Given an orthonormal basis~$\omega_1,\ldots, \omega_{2\mu}$
of~$\bigwedge^1 T$, and real numbers~$\lambda_1, \ldots, \lambda_\mu$,
the form
\begin{equation}
\label{25}
\form = \sum_{j=1}^\mu \lambda_j \left( \omega_{2j-1} \wedge
\omega_{2j} \right)
\end{equation}
has comass~$\|\form \|=\max_j |\lambda_j|$.
\end{lemma}

\begin{proof}
We can assume without loss of generality that each~$\lambda_j$ is
nonnegative.  This can be attained in one of two ways.  One can
permute the coordinates, by applying the transposition flipping
$\omega_{2j-1}$ and~$\omega_{2j}$.  Alternatively, one can replace,
say,~$\omega_{2j}$ by~$-\omega_{2j}$.

Next, consider the hermitian inner product~$H_f$ obtained by
polarizing the quadratic form
\[
\sum_j \left(\lambda_j^{1/2} \omega_{2j}\right)^2+
\left(\lambda_j^{1/2} \omega_{2j+1}\right)^2.
\]
Let~$\zeta= v\wedge w$ be an orthonormal pair such that~$||\form||=
\form(\zeta)$.  As in~\eqref{21}, we have
\[
f(\zeta)=-iH_f(\zeta) =H_f(iv,w) \leq \left( \max_j \lambda_j \right)
(iv)\cdot w \leq \max_j \lambda_j,
\]
proving the lemma.
\end{proof}

\section{Gromov's inequality for complex projective space}
\label{gro}

First we recall the definition of the stable norm in the real
$k$-homology of an~$n$-dimensional polyhedron~$X$ with a piecewise
Riemannian metric, following \cite{BK1, BK2}.

\begin{definition}
\label{D1.2}
The stable norm~$\parallelh$ of~$h\in H_k(X,\RR)$ is the infimum of
the volumes
\begin{equation}
\vol_k(c)=\Sigma_i |r_i| \vol_k(\sigma_i)
\end{equation}
over all real Lipschitz
cycles~$c=\Sigma_i r_i \sigma_i$ representing~$h$.
\end{definition}

Note that~$\|\;\|$ is indeed a norm, \cf~\cite{Fe2} and
\cite[4.C]{Gr3}.

We denote by~$H_k(X,\Z)_{\rr}$ the image of~$H_k(X,\Z)$ in
$H_k(X,\RR)$ and by~$h_{\rr}$ the image of~$h\in H_k(X,\Z)$ in
$H_k(X,\RR)$.  Recall that~$H_k(X,\Z)_{\rr}$ is a lattice in
$H_k(X,\RR)$.  Obviously
\begin{equation}
\label{D(1.3)}
\parallelhrr \le \vol_k(h)
\end{equation}
for all~$h\in H_k(X,\Z)$, where~$\vol_k(h)$ is the infimum of volumes
of all integral~$k$-cycles representing~$h$.  Moreover, one
has~$\parallelhrr = \vol_n(h)$ if~$h\in H_n(X,\Z)$.  H.~Federer
\cite[4.10, 5.8, 5.10]{Fe2} (see also \cite[4.18 and~4.35]{Gr3})
investigated the relations between~$\parallelhrr$ and~$\vol_k(h)$ and
proved the following.
\begin{proposition}
If~$h\in H_k(X,\Z)$,~$1\le k < n$, then
\begin{equation}
\label{1013b}
\parallelhrr = \lim\limits_{i\rightarrow\infty} \frac{1}{i} \vol_k (i
h).
\end{equation}
\end{proposition}

Equation \eqref{1013b} is the origin of the term {\em stable
norm\/}\index{stable norm} for~$\parallelleer$.  Recall that the
stable~$k$-systole of a metric~$(X,\gmetric)$ is defined by setting
\begin{equation}
\label{11zbis}
\stsys_k(\gmetric)=\lambda_1\left( H_k(X, \Z)_\R^{\phantom{I}}, \|\;\|
\right),
\end{equation}
\cf \eqref{12b} and \eqref{11z}.  Let us now return to systolic
inequalities on projective spaces.

\begin{theorem}[M.~Gromov]
\label{1211}
Every Riemannian metric~$\gmetric$ on complex projective space~$\CP^n$
satisfies the inequality
\[
\stsys_2(\gmetric)^n \leq n! \vol_{2n}(\gmetric);
\]
equality holds for the Fubini-Study metric on~$\CP^n$.
\end{theorem}
\begin{proof}
Following Gromov's notation in \cite[Theorem 4.36]{Gr3}, we let
\begin{equation}
\label{31}
\alpha\in H_2(\CP^n;\Z)=\Z
\end{equation}
be the positive generator in homology, and let
\[
\omega\in H^{2}(\CP^n;\Z)=\Z
\]
be the dual generator in cohomology.  Then the cup power~$\omega^n$ is
a generator of~$H^{2n}(\CP^n;\Z)=\Z$.  Let~$\eta \in \omega$ be a
closed differential 2-form.  Since wedge product~$\wedge$
in~$\Omega^*(X)$ descends to cup product in~$H^*(X)$, we have
\begin{equation}
\label{P61b}
1= \int_{\CP^{n}} \eta^{\wedge n} .
\end{equation}
Now let~$\gmetric$ be a metric on~$\CP^n$.

The comass norm of a differential~$k$-form is, by definition, the
supremum of the pointwise comass norms, \cf Definition~\ref{34}.  Then
by the Wirtinger inequality and Corollary~\ref{1214}, we obtain
\begin{equation}
\label{P61d}
\begin{aligned}
1 & \leq \int_{\CP^n} \| \eta^{\wedge n}\|\; d\!\vol \\ & \leq n!
\left( \| \eta \|_\infty \right)^n \vol_{2n}(\CP^n,\gmetric)
\end{aligned}
\end{equation}
where~$\| \; \|_\infty$ is the comass norm\index{comass norm} on forms
(see \cite[Remark~4.37]{Gr3} for a discussion of the constant in the
context of the Wirtinger inequality\index{Wirtinger inequality}).  The
infimum of \eqref{P61d} over all~$\eta\in \omega$ gives
\begin{equation}
\label{P61c}
1 \leq n! \left( \| \omega \|^* \right)^n \vol_{2n} \left( \CP^n,
\gmetric \right),
\end{equation}
where~$\|\;\|^*$ is the comass norm in cohomology.  Denote by~$\|\;\|$
the stable norm\index{stable norm} in homology.  Recall that the
normed lattices~$(H_2(M;\Z), \|\;\|)$ and~$(H^2(M;\Z), \|\;\|^*)$ are
dual to each other \cite{Fe1}.  Therefore the class~$\alpha$
of~\eqref{31} satisfies
\[
\|\alpha \| = \frac{1}{\|\omega \|^*},
\]
and hence
\begin{equation}
\label{P61}
\stsys_2(\gmetric)^n = \| \alpha \|^n \leq n! \vol_{2n}(\gmetric) .
\end{equation}
Equality is attained by the two-point homogeneous Fubini-Study metric,
since the standard~$\CP^1 \subset \CP^n$ is calibrated by the
Fubini-Study Kahler 2-form, which satisfies equality in the Wirtinger
inequality at every point.
\end{proof}

\begin{example}
\label{1222}
Every metric~$\gmetric$ on the complex projective plane satisfies the
optimal inequality
\[
\stsys_2(\CP^2,\gmetric)^2 \leq 2 \vol_{4}(\CP^2,\gmetric).
\]
\end{example}

This example generalizes to the manifold obtained as the connected sum
of a finite number of copies of~$\CP^2$ as follows.

\begin{proposition}
Every Riemannian~$n\cptwo$ satisfies the inequality
\begin{equation}
\label{N42}
\stsys_2\left(n\cptwo \right)^2\leq 2 \vol_{4}\left( n\cptwo \right).
\end{equation}
\end{proposition}

\begin{proof}
We define two varieties of conformal~$2$-systole of a manifold~$M$ as
follows.  The Euclidean norm~$|\;|$ and the comass norm~$\|\;\|$ on
(linear)~$2$-forms define, by integration, a pair of~$L^2$ norms on
$\Omega^2(M)$.  Minimizing over representatives of a cohomology class,
we obtain a pair of norms in de Rham cohomology.  The dual norms in
homology will be denoted respectively~$|\;|_2$ and~$\|\;\|_2$, \cf
\cite[p.~122, 130]{SGT}.  We let
\[
\Confsys_2= \lambda_1(H_2(M;\Z), \|\;\|_2)
\]
and
\[
\confsys_2= \lambda_1(H_2(M;\Z), |\;|_2).
\]
Since every top dimensional form is simple (decomposable), by
Corollary~\ref{1214} we have an inequality
\begin{equation}
\label{310}
|x|^2 \leq \Wirt_2 \|x\|^2
\end{equation}
where~$\Wirt_2=2$, between the pointwise Euclidean norm and the
pointwise comass, for all~$x\in \bigwedge^2(n\CP^2)$.  It follows
that, dually, we have
\begin{equation}
\label{311}
        \Confsys_2^2 \leq 2 \confsys_2^2.
\end{equation}
For a metric of unit volume we have
\begin{equation}
\label{312}
\stsys_k \leq \Confsys_k.
\end{equation}
Combining \eqref{311} and \eqref{312}, we obtain
\[
        \stsys_2^2(\gmetric) \leq 2 \confsys_2^2(\gmetric)
        \vol_4(\gmetric) .
\]
Recall that the intersection form of~$n\CP^2$ is given by the identity
matrix.  Every metric~$\gmetric$ on a connected sum~$n\CP^2$ satisfies
the identity~$\confsys_2(\gmetric)=1$ because of the identification of
the~$L^2$ norm and the intersection form.  We thus reprove Gromov's
optimal inequality
\[
        \stsys_2^2 \leq 2 \vol_4,
\]
but now it is valid for the connected sum of~$n$ copies of~$\CP^2$.
\end{proof}

In fact, the inequality can be stated in terms of the last successive
minimum~$\lambda_n$ of the integer lattice in homology with respect to
the stable norm~$\|\;\|$.

\begin{corollary}
The last successive minimum~$\lambda_n$ satisfies the inequality
\[
\lambda_n \left( H_2(n\CP^2,\Z), ||\;|| \right) ^2 \leq 2
\vol_4(n\CP^2)
\]
\end{corollary}
The proof is the same as before.  This inequality is in fact optimal
for all~$n$, though equality may not be attained.

\begin{question}
What is the asymptotic behavior for the stable systole of~$n\CP^2$
when~$n\to \infty$?  Can the constant in \eqref{N42} be replaced by a
function which tends to zero as~$n\to \infty$?
\end{question}

\section{Symmetric metric of~$\HP^2$ and Kraines~$4$-form}
\label{1215}

The quaternionic projective plane~$\HP^2$ has volume~$\vol_8(\HP^2) =
\tfrac{\pi^4}{5!}$ for the symmetric metric with sectional curvature
$1\leq K \leq 4$, while for the projective line with~$K\equiv 4$ we
have~$\vol_4(\HH \PP^1)= \tfrac{\pi^2}{3!}$, \cf
\cite[formula~(3.10)]{Be1}.  Since the projective line is volume
minimizing in its real homology class, we obtain~$\stsys_4(\HH \PP^2)=
\tfrac{\pi^2}{3!}$, as well, resulting in a systolic ratio
\begin{equation}
\label{1221}
\frac {\stsys_4(\HP^2)^2} {\vol_8(\HP^2)} = \tfrac{10}{3}
\end{equation}
for the symmetric metric.

In more detail, we endow~$\HP^n$ with the natural metric as the base
space of the Riemannian submersion from the unit
sphere
\[
S^{4n+3}\subset \HH^{n+1} .
\]
A projective line~$\HP^1 \subset \HP^n$ is a round~$4$-sphere of
(Riemannian) diameter~$\tfrac{\pi}{2}$ and sectional curvature~$+4$,
attaining the maximum of sectional curvatures of~$\HP^n$.  The
extension of scalars from~$\R$ to~$\HH$ gives rise to an
inclusion~$\R^3\hookrightarrow \HH^3$, and thus an
inclusion~$\RP^2\hookrightarrow \HP^2$.  Then~$\RP^2\subset \HP^2$ is
a totally geodesic submanifold of diameter~$\tfrac{\pi}{2}$ and
Gaussian curvature~$+1$, attaining the minimum of the sectional
curvatures of~$\HP^2$, \cf \cite[p.~73]{CE}.

The following proposition was essentially proved by V.~Kraines
\cite{Kr} and M.~Berger \cite{Be1}.  The invariant 4-form was briefly
discussed in \cite[p.~152]{HL}.

\begin{proposition}
There is a parallel~$4$-form~$\kraines\in \Omega^4(\HP^2)$
representing a generator of~$H^4(\HP^2,\Z)=\Z$, with
\begin{equation}
\label{1232}
|\kappa_\HP^2| =\tfrac{10}{3} \|\kraines\|^2
\end{equation}
and
\begin{equation}
\label{1233}
|\kraines|^2 = \tfrac{10}{3}\|\kraines\|^2,
\end{equation}
where~$|\;|$ and~$\|\;\|$ are, respectively, the Euclidean norm and
the comass of the unit volume symmetric metric on~$\HP^2$.
\end{proposition}

\begin{proof}
The parallel differential~$4$-form~$\kraines$ is obtained from an
$Sp(2)$-invariant alternating~$4$-form on a tangent space at a point,
by propagating it by parallel translation to all points of~$\HP^2$.
The fact that parallel translation produces a well defined global
$4$-form results from the~$Sp(2)$ invariance of the alternating form.

In more detail, consider the quaternionic vector
space~$\HH^n=\R^{4n}$.  Each of the three quaternions~$i$,~$j$,
and~$k$ defines a complex structure on~$\HH^n$, \ie an
identification~$\HH^n \simeq \C^{2n}$.  The imaginary part of the
associated Hermitian inner product on~$\C^{2n}$ is the standard
symplectic exterior~$2$-form, \cf Example~\ref{1211z}.
Let~$\omega_i$,~$\omega_j$, and~$\omega_k$ be the triple of~$2$-forms
on~$\HH^n$ defined by the three complex structures.  We consider their
wedge squares~$\omega_i^2$,~$\omega_j^2$, and~$\omega_k^2$.  We define
an exterior~$4$-form~$\krainesl_n$, first written down explicitly by
V.~Kraines~\cite{Kr}, by setting
\begin{equation}
\label{krai}
\krainesl_n= \tfrac{1}{6} \left( \omega_i^2+\omega_j^2+\omega_k^2
\right).
\end{equation}
The coefficient~$\tfrac{1}{6}$ normalizes the form to unit comass, \cf
Lemma~\ref{1213}.  The form~$\krainesl_n$ is invariant under
transformations in~$Sp(n)\times Sp(1)$ \cite[Theorem~1.9]{Kr} and thus
defines a parallel differential~$4$-form in~$\Omega^4\HP^n$, which is
furthermore closed.  We normalize the differential form in such a way
as to represent a generator of integral cohomology, and denote the
resulting form~$\kraines$, so that~$[\kraines]\in H^4(\HP^n,\Z)_\R
\simeq \Z$ is a generator.

In the case~$n=2$, explicit formulas appear in \eqref{81b} and
\eqref{cartanbasis}.  Here~$\omega_i$ is the sum of~$4$ monomial
terms, while~$\omega_i^2$ is twice the sum of~$6$ such terms.

The form~$3 \krainesl_2$ on~$\HH^2$ decomposes into a sum of~$18$
simple~$4$-forms, \ie monomials in the~$8$ coordinates.  The~$18$
monomials are not all distinct.  Two of them, denoted~$m_0$ and its
Hodge star~$* m_0$, occur with multiplicity~$3$.  Thus, we obtain a
decomposition as a linear combination of seven selfdual pairs
\begin{equation}
\label{seven}
3\krainesl_2= 3(m_0 + * m_0) + \sum_{\ell=1}^6 (m_\ell + *
m_\ell),
\end{equation}
where~$*$ is the Hodge star operator.  In Section~\ref{eight}, the
explicit formulas for the three~$2$-forms will be used to write down
the Cayley~$4$-form~$\cay$.

Similarly to \eqref{P61d}, we can write
\begin{equation}
\label{49}
\begin{aligned}
1 & = \int_{\HP^2} \left| \kraines^{\wedge 2} \right|\; d\!\vol \\ & =
\tfrac{10}{3} \left( \| \kraines \|_\infty \right)^2 \vol_{8}(\HP^2),
\end{aligned}
\end{equation}
thereby reproving \eqref{1221} by the duality of comass and stable
norm.

\begin{lemma}
The Kraines form~$\kappa_2$ of \eqref{krai} has unit comass:~$\left\|
\krainesl_2 \right\| = 1$.
\end{lemma}

This was proved in \cite{Be1,DHM}.  Meanwhile, from \eqref{seven} we
have
\[
\left( 3\krainesl_2 \right)^2= 2 \left( 9 \vol + 6 \vol
\right),
\]
where~$\vol=e_1 \wedge e_2 \wedge \cdots \wedge e_8$ is the volume
form of~$\HH^2=\R^8$.  Hence
\[
\left| \left( 3\krainesl_2 \right)^2 \right|= 2\cdot 15=30,
\]
proving identity~\eqref{1232}.  Meanwhile,~$\left| 3
\krainesl_2 \right|^2 = 9 + 9 + 12 =30$, proving
identity~\eqref{1233}.
\end{proof}

\begin{remark}
There is a misprint in the calculation of the systolic constants in
\cite[Theorem~6.3]{Be1}, as is evident from
\cite[formula~(6.14)]{Be1}.  Namely, in the last line on page
\cite[p.~12]{Be1}, the formula for the coefficient~$s_{4,b}$ lacks the
exponent~$b$ over the constant~$6$ appearing in the numerator.  The
formula should be
\[
s_{4,b}= \frac{6^b}{(2b+1)!}.
\]
\end{remark}

\section{Generalized Wirtinger inequalities}
\label{e7}

\begin{definition}
The {\em Wirtinger constant\/}~$\W_n$ of~$\R^{2n}$ is the maximal
ratio~$\tfrac{|\omega^2|}{\|\omega\|^2}$ over all~$n$-forms~$\omega\in
\Lambda^n \R^{2n}$.  The {\em modified Wirtinger constant\/}~$\W'_n$
is the maximal ratio~$\tfrac{|\omega|^2}{\|\omega\|^2}$ over
$n$-forms~$\omega$ on~$\R^{2n}$.
\end{definition}

The calculation of~$\W_n$ can thus be thought of as a generalisation
of the Wirtinger inequality of Section~\ref{fed}.

In Section~\ref{nine}, we will deal in detail with the special case of
self-dual~$4$-forms in the context of the Lie algebra~$E_7$.  We
therefore gather here some elementary material pertaining to this
case.

\begin{definition}
Let~$n$ be even.  Let~$\W_{sd}$ be the maximal
ratio~$\tfrac{|\omega^2|}{\|\omega\|^2}$ over all selfdual~$n$-forms
on~$\R^{2n}$.
\end{definition}

\begin{lemma}
\label{121}
One has~$\W_n =\W_{sd} \leq \W_n'$ if~$n$ is even.
\end{lemma}

\begin{proof}
In general for a skew-form~$\omega$ it may occur that
$|\omega^2|>|\omega|^2$.  This does not occur when~$\omega$ is
middle-dimensional.  If~$\omega$ is a middle-dimensional form, then
\begin{equation}
\label{1234}
\| \omega ^2 \| = | \omega^2 | = \langle \omega, *\omega \rangle \leq
|\omega| \; |\!*\omega| = |\omega|^2,
\end{equation}
proving that~$\W_n \leq \W'_n$.

Let~$\eta$ be a form with nonnegative wedge-square (if it is negative,
reverse the orientation of the ambient vector space~$\R^{2n}$ to make
the square non-negative, without affecting the values of the relevant
ratios).  If~$n$ is even, the Hodge star is an involution.
Let~$\eta=\eta_+ + \eta_-$ be the decomposition into selfdual and
anti-selfdual parts under Hodge~$*$.  Then
\begin{equation}
\label{P61dz}
\begin{aligned}
\eta^{2} & = \left( \eta_+ + \eta_- \right)^{2} \\ & = \eta_+^2 +
\eta_-^2
\end{aligned}
\end{equation}
Thus
\begin{equation}
|\eta^{2}| = | \eta_+^2 |- | \eta_-^2 | \leq | \eta_+^2 | .
\end{equation}
Meanwhile,
\[
\|\eta_+ \| = \tfrac{1}{2} \left( \|\eta + *\eta \| \right) \leq
\tfrac{1}{2} \left( \|\eta \| + \| *\eta \|\right) = \|\eta\|
\]
by the triangle inequality.  Thus,~$\| \eta_+ \| \leq \| \eta \|$ and
we therefore conclude that
\[
\frac{|\eta^2|}{\|\eta\|^2} \leq \frac{|\eta_+^2|}{\|\eta_+\|^2}
\leq \W_{sd},
\]
proving that~$\W_n=\W_{sd}$.
\end{proof}

\begin{proposition}
\label{upperbound}
Let~$X$ be an orientable, closed manifold of dimension~$2n$, with
$b_n(X)=1$.  Then
\[
\SR_n(X) \leq \W_n.
\]
\end{proposition}
\begin{proof}

By Poincar\'e duality, the fundamental cohomology class in the group
$H^{2n}(X;\Z) \simeq \Z$ is the cup square of a generator of the
cohomology group~$H^n(X;\Z)_\R \simeq \Z$.  The inequality is now
immediate by applying the method of proof of \eqref{P61d}.
\end{proof}

Recall that the cohomology ring for~$\CP^n$ is polynomial on a single
$2$-dimensional generator, truncated at the fundamental class.  The
cohomology ring for~$\HP^n$ is the polynomial ring on a single
$4$-dimensional generator, similarly truncated.  Thus the middle
dimensional Betti number is 1 if~$n$ is even and 0 if~$n$ is odd.

\begin{corollary}
\label{55}
Let~$n\in \N$.  We have the following bounds for the
middle-dimensional stable systolic ratio:
\[
\begin{aligned}
\SR_{4n}(\HP^{2n}) & \leq \W_{4n} \\ \SR_{2n}(\CP^{2n}) & \leq \W_{2n}
\\ \SR_{8}(M^{16}) & \leq \W_{8}
\end{aligned}
\]
where~$M^{16}$ is the Cayley projective plane.
\end{corollary}

\begin{remark}
The systolic ratio of the symmetric metric of~$\CP^4$ is~$6$, while by
Proposition~\ref{45} we have~$\W_4=14 > 6$, so that Corollary~\ref{55}
gives a weaker upper bound of~$14$ for the optimal systolic ratio
of~$\CP^4$.  Thus it is in principle impossible to calculate the
optimal systolic ratio for either~$\HP^2$ or~$\CP^4$ by any direct
generalisation of Gromov's calculation \eqref{P61d}.
\end{remark}

The detailed calculation of the Wirtinger constant~$\W_4$ appears in
Section~\ref{nine}.

\section{$BG$ spaces and a homotopy equivalence}
\label{beegee}

Systolically interesting metrics can be constructed as pullbacks by
homotopy equivalences.  A particularly useful one is described below.

\begin{proposition}
The complex projective~$2n$-space~$\CP^{2n}$ admits a degree~$1$ map
to the quaternionic projective space~$\HP^n$.
\end{proposition}

\begin{proof}
Such a map can be defined in coordinates by including~$\C^{2n+1}$
in~$\C^{2n+2}$ as a hyperplane, identifying~$\C^{2n+2}$
with~$\HH^{n+1}$, and passing to the appropriate quotients.  To verify
the assertion concerning the degree in a conceptual fashion, we
proceed as follows.  We imbed~$\CP^{2n}$ as the~$(4n)$-skeleton
of~$\CP^\infty$.  The latter is a model for the classifying
space~$BS^1$ of the circle.  Similarly, we have
\[
\HP^n = \left( \HP^\infty \right)^{(4n)} \subset \HP^\infty \simeq
BS^3,
\]
where~$S^3$ is identified with the unit quaternions.  Namely,~$BG$ can
be characterized as the quotient of a contractible space~$S$ by a
free~$G$ action.  But~$\HP^\infty$ is such a quotient for~$S=S^\infty$
and ~$G = S^3$.  The inclusion of~$S^1$ as a subgroup of~$S^3$ defines
a map~$\CP^\infty \to \HP^\infty$.  The composed map~$\CP^{2n}
\hookrightarrow \CP^\infty \to \HP^\infty$ is compressed, using the
cellular approximation theorem, to the~$(4n)$-skeleton.  In matrix
terms, an element~$u\in S^1$ goes to the element
\begin{equation}
\label{61}
\left[\begin{matrix} u & 0 \\ 0 & u^{-1}
\end{matrix} \right] \in SU(2)=S^3.
\end{equation}

The induced map on cohomology is computed for the infinite dimensional
spaces, and then restricted to the~$(4n)$-skeleta.  By
Proposition~\ref{61b}, the cohomology of~$BS^3$ is~$\Z[c_2]$, \ie a
polynomial algebra on a~$4$-dimensional generator~$c_2$, given by the
second Chern class.  Thus, to compute the induced homomorphism
on~$H^4$, we need to compute~$c_2$ of the sum of the tautological line
bundle~$L$ on~$\CP^\infty$ and its inverse, \cf \eqref{61}.  By the
sum formula, it is
\[
-c_1(L)^2,
\]
but this is a generator of~$H^4(\CP^\infty)$.  In other words, the map
\[
H^4(BS^3) \to H^4(BS^1)
\]
is an isomorphism.  From the structure of the cohomology algebra, we
see that the same is true for the induced homomorphism in~$H^{4n}$.
The inclusions of the~$(4n)$-skeleta of these~$BG$ spaces are
isomorphisms on cohomology~$H^{4n}$, as well, in view of the absence
of odd dimensional cells.  Hence the conclusion follows for these
finite-dimensional projective spaces.
\end{proof}

The lower bound of Theorem~\ref{11} for the optimal systolic ratio
of~$\HP^2$ follows from the two propositions below.

\begin{proposition}
\label{61b}
We have~$H^*(B S^3) = \Z[v]$, where the element~$v$
is~$4$-dimensional.  Meanwhile,~$H^*(BS^1) =\Z[c]$, where~$c$
is~$2$-dimensional.  Here~$i^*(v) = -c^2$ (with usual choices for
basis),~$S^3 = SU(2)$, and~$v$ is the second Chern class.
\end{proposition}

Now restrict attention to the~$4n$-skeleta of these spaces.  We obtain
a map
\begin{equation}
\label{51}
\CP^{2n} \to \HP^n
\end{equation}
which is degree one (from the cohomology algebra).

\begin{proposition}
\label{62}
There exists a map~$\HP^n \to \CP^{2n} \cup e^3 \cup e^7 \cup \ldots
\cup e^{4n-1}$ defining a homotopy equivalence.
\end{proposition}

\begin{proof}
Coning off a copy of~$\CP^1\subset \CP^{2n}$, we note that the map
\eqref{51} factors through the CW complex~$\CP^{2n}\cup e^3$.

The map~$\CP^4\cup e^3 \to \HP^2$ is an isomorphism on homology
through dimension~$5$, and a surjection in dimension~$6$.  We consider
the pair
\[
(\HP^2, \CP^4\cup e^3).
\]
Its homology vanishes through dimension~$6$ by the exact sequence of a
pair.  The relative group~$H_7(\HP^2, \CP^4\cup e^3)$ is mapped by the
boundary map to~$H_6(\CP^4\cup e^3) = \Z$, generated by an
element~$h\in H_6(\CP^4\cup e^3)$.  We therefore obtain an isomorphism
\[
\alpha: H_6(\CP^4\cup e^3) \to H_7(\HP^2, \CP^4\cup e^3),
\]
\cf Figure~\ref{Cfig1}.

\renewcommand{\arraystretch}{1.3}
\begin{figure}
\[
\xymatrix { H_6(\CP^4\cup e^3) \ar[d]_\alpha \ar[r] & \pi_6(\CP^4\cup
e^3) \\ H_7(\HP^2,\CP^4\cup e^3) \ar[r]^{\beta^{-1}} &
\pi_7(\HP^2,\CP^4\cup e^3) \ar[u]_\gamma }
\]
\caption{\textsf{Commutation of boundary and Hurewicz homomorphisms}}
\label{Cfig1}
\end{figure}
\renewcommand{\arraystretch}{1}

Both spaces are simply connected and the pair is 6-connected as a
pair.  Applying the relative Hurewicz theorem, we obtain an
isomorphism
\[
\beta: \pi_7(\HP^2, \CP^4\cup e^3) \to H_7(\HP^2, \CP^4\cup e^3).
\]
Applying the boundary homomorphism
\[
\gamma: \pi_7(\HP^2, \CP^4\cup e^3) \to \pi_6(\CP^4\cup e^3),
\]
we obtain an element
\begin{equation}
\label{63}
h'= \gamma \circ \beta^{-1} \circ \alpha(h)\in \pi_6(\CP^4\cup e^3)
\end{equation}
which generates~$H_6$ and is mapped to~$0\in \pi_6(\HP^2)$.

We now attach a~$7$-cell to the complex~$\CP^4\cup e^3$ using the
element~$h'$ of \eqref{63}.  We obtain a new CW complex
\[
X = \left( \CP^4 \cup e^3 \right)\cup _{h'} e^7,
\]
and a map~$X \to \HP^2$, by choosing a nullhomotopy of the composite
map to~$\HP^2$.  The new map is an isomorphism on all homology.  Since
both spaces are simply connected, the map is a homotopy equivalence.
Reversing the arrow, we obtain a homotopy equivalence from~$\HP^2$ to
the union of~$\CP^4$ with cells of dimension~$3$ and~$7$.  A similar
argument, applied inductively, establishes the general case.
\end{proof}

\forget

\begin{proposition}
There is a hands-on construction of a map
\[
    \HP^k \to \CP^{2k} \cup e^3 \cup e^7 \cup \ldots ,
\]
which splits the twistor bundle
\[
    p: \CP^{2k+1} \to \HP^k,
\]
restricted to the~$2k$-skeleton.
\end{proposition}

\begin{proof}
The construction is by induction on the skeleta as follows.  Note that
the twistor bundle respects the skeleta, \ie the restriction of the
projection~$p$ to~$\CP^{2m+1}$, where~$m < k$, will be the
corresponding twistor bundle over~$\HP^m \subset \HP^k$.  We therefore
proceed as follows.

Step 1.  We attach a~$3$-disk to~$\CP^1 \subset \CP^2$.  This allows
us to construct a section of the twistor bundle over~$\HP^1 = S^4$.
The complex thus obtained is denoted~$X_4$.  The complex~$X_4$ is
homotopy equivalent to~$\HP^1$.

Step 2.
Consider the space
\[
    X_6 = X_4 \cup \CP^3,
\]
where the union is along a common~$\CP^2$.  Obviously, projection~$p$
of the twistor bundle extends from~$\CP^3$ to~$X_6$, and one can view
$X_6$ as~$X_4 \cup e^6$, where~$e^6$ is the cell of maximal dimension
in~$\CP^3$.

Now note that this cell is attached along a homotopically trivial map
$S^5 \to X_4$.  Indeed, the image~$p(e^6) \subset \HP^1 \simeq X_4$
provides the required homotopy.

Thus, we can attach a~$7$-cell~$e^7$ to~$X_6$ in such a way that the
resulting complex~$X_7$ admits a retraction to its subcomplex~$X_4$.
Since we have an inclusion~$\CP^3 \subset X_7$, we can naturally
define the complex
\[
X_8 = X_7 \cup \CP^4,
\]
where the union is along a common~$\CP^3$.  In view of the deformation
retraction of~$X_7$ to~$X_4$, the map~$p: \CP^4 \to \HP^2$ extends to
the entire complex~$X_8$.  It is easy to see that this map induces
isomorphism in all homology.  Since both spaces are simply connected,
we obtain a homotopy equivalence.

Step 3. This argument can be easily extended by induction to all
dimensions, as follows.  If the space
\[
X_{4k} = \CP^{2k} \cup e^3 \cup e^7 \ldots \cup e^{4k-1}
\]
has already been constructed, then the~$(4k+2)$-dimensional cell in
$X_{4k+2} = X_{4k} \cup \CP^{2k+1}$ is always attached along a
homotopically trivial map, since we have a map~$p : X_{4k+2} \to \HP^k
\simeq X_{4k}$.  It is "ubivaetsia" by attaching a cell~$e^{4k+3}$ in
the same way is in Step 2 above.
\end{proof}

\forgotten

\section{Lower bound for quaternionic projective space}
\label{sevenb}

In this section, we apply the homotopy equivalence constructed in
Section~\ref{beegee}, so as to obtain systolically interesting
metrics.

\begin{proposition}
\label{71}
One can homotope the map of Proposition~\ref{62} to a simplicial map,
and choose a point in a cell of maximal dimension in
\begin{equation}
\label{71v}
\CP^{2n} \subset \CP^{2n} \cup e^3 \cup \ldots \cup e^{4n-1}
\end{equation}
with a unique inverse image.
\end{proposition}

\begin{proof}
To fix ideas, consider the case~$n=2$.  The inverse image of a little
ball around such a point is a union of balls mapping the obvious way
to the ball in~$\CP^4 \cup e^3 \cup e^7$.  We need to cancel balls
occurring with opposite signs.  Take an arc connecting the boundaries
of two such balls where the end points are the same point of the
sphere.  Apply homotopy extension to make the map constant on a
neighborhood of this arc ($\pi_1$ of the target is 0).  Then the union
of these balls and fat arc is a bigger ball and we have a
nullhomotopic map to the sphere on the boundary.  We can homotope the
map to the disc relative to the boundary to now lie in the sphere.
\end{proof}

\begin{corollary}
\label{53}
The optimal middle dimensional stable systolic ratio of~$\HP^{2n}$
equals that of~$\CP^{4n}$.
\end{corollary}

\begin{proof}
We first prove the inequality~$\SR_{4n}(\CP^{4n}) \geq
\SR_{4n}(\HP^{2n})$.  We exploit the degree one map~\eqref{51}.
Recall that a map is called {\em monotone\/} if the preimage of every
connected set is connected.  By the work of A.~Wright \cite{Wr}, the
map~\eqref{51} can be homotoped to a simplicial monotone map.  In
particular, the preimage of every top-dimensional simplex is a single
top-dimensional simplex.  Thus the pull-back ``metric'' has the same
total volume as the metric of the target.  Pulling back metrics
from~$\HP^{2n}$ to~$\CP^{4n}$ by the monotone simplicial map completes
the proof in this direction.

Let us prove the opposite inequality.  To fix ideas, we let~$n=1$.  We
need to show that~$\SR_4(\CP^4) \leq \SR_4(\HP^2)$.  Once the map

\begin{equation}
\label{eff}
f: \HP^2 \to \CP^4 \cup e^3 \cup e^7
\end{equation}
is one-to-one on an~$8$-simplex
\[
\Delta\subset \CP^4 \cup e^3 \cup e^7
\]
of the target (by Proposition~\ref{71}), we argue as follows.  The
images of the attaching maps of~$e^3$ and~$e^7$ may be assumed to lie
in a hyperplane~$\CP^3\subset \CP^4$.  Take a self-diffeomorphism
\begin{equation}
\label{fi}
\phi: \CP^4 \to \CP^4
\end{equation}
preserving the hyperplane, and sending the 8-simplex~$\Delta$ to the
complement of a thin neighborhood of the hyperplane, so that most of
the volume of the symmetric metric of~$\CP^4$ is contained in the
image of~$\Delta$.

Now pull back the metric of the target by the composition~$\phi \circ
f$ of the maps \eqref{eff} and \eqref{fi}.  The resulting ``metric''
on~$\HP^2$ is degenerate on certain simplices.  The metric can be
inflated slightly to make the quadratic form nondegenerate everywhere,
without affecting the total volume significantly.  The proof is
completed by the following proposition.
\end{proof}

\begin{proposition}
Fix any background metric on~$\CP^{4n}$, e.g. the Fubini-Study.  Then
the metric can be extended to the~$3$-cell, the~$7$-cell,~$\ldots,$
the~$(8n-1)$-cell, as in \eqref{71v}, in such a way as to decrease the
stable systole by an arbitrarily small amount.
\end{proposition}

\begin{proof}
We work in the category of simplicial polyhedra~$X$, \cf \cite{Ba06}.
Here volumes and systoles are defined, as usual, simplex by simplex.
When attaching a cell along its boundary, the attaching map is always
assumed to be simplicial, so that all systolic notions are defined on
the new space, as well.

The metric on the attached cells needs to be chosen in such a way as
to contain a long cylinder capped off by a hemisphere.

\forget
Since a minimizing integral 4-current would avoid a 3-cell, it is the
7-cell that poses a potential problem.  The cylinder in this case is a
product of an interval by a sphere~$S^6$.  We choose a metric on the
cylinder large enough so as to satisfy two conditions.  On the one
hand, the complement of the cap admits a distance-decreasing map back
to~$\CP^4$.

On the other hand, we choose the metric on the cylinder in such a way
that the cylinder has sufficiently small curvature, so that the
monotonicity formula, applied to a metric ball contained in the
cylinder, will give an area lower bound larger than the stable
4-systole of~$\CP^4$.  This ensures that any minimizing integral
current representing a generator of~$H_4$ will not enter the 7-cell.

However, this is not quite enough, since, by the stabilisation
formula, the calculation of the stable systole of the CW complex
involves the volumes of minimizing currents from multiple classes, as
well.
\forgotten

To make sure the attachment of a cell~$e^p$ does not significantly
decrease the stable systole, we argue as follows.

To fix ideas, let~$n=1$.  Normalize~$X$ to unit stable 4-systole.
Let~$W= X \cup e^p$, and consider a metric on~$e^p$ which includes a
cylinder of length~$L>>0$, based on a sphere~$S^{p-1}$, of radius~$R$
chosen in such a way that the attaching map~$\partial e^p \to X$ is
distance-decreasing.  Here~$R$ is fixed throughout the argument (and
in particular is independent of~$L$).

Now consider an~$n$-fold multiple of the generator~$g\in H_4(W)$, well
approximating the stable norm in the sense of \eqref{1013b}.  Consider
a simplicial~$4$-cycle~$M$ with integral coefficients, in the
class~$ng\in H_4(W)$.  We are looking for a lower bound for the stable
norm~$\|g\|$ in~$W$.  Here we have to deal with the possibility that
the 4-cycle~$M$ might ``spill'' into the cell~$e^p$.  Applying the
coarea inequality~$\vol_4(M) \geq \int_0^L \vol_3(M_t) dt$ along the
cylinder, we obtain a 3-dimensional section~$S=M_{t_0}$ of~$M$
of~$3$-volume at most
\begin{equation}
\label{73}
\vol_3(S) = \frac{n \|g\|}{L},
\end{equation}
\ie as small as one wishes compared to the~$4$-volume of~$M$ itself.
Here~$M$ decomposes along~$S$ as the union
\[
        M = M_+ \cup M_-
\]
where~$M_+$ admits a distance decreasing projection to the
polyhedron~$X$, while~$M_-$ is entirely contained in~$e^p$.  For
any~$4$-chain~$C\subset S^{p-1}$ filling~$S$, the new~$4$-cycle
\[
        M' = M_+ \cup C
\]
represents the same homology class~$ng\in H_4(W)$, since the
difference~$4$-cycle~$M - M'$ is contained in a~$p$-ball whose
homology is trivial.  Now we apply the linear (without the
exponent~$\frac{n+1}{n}$) isoperimetric inequality in~$S^{p-1}$.  This
allows us to fill the section~$S=\partial M_+$ by a
suitable~$4$-chain~$C\subset S^{p-1}$ of volume at most
\[
\vol_4(C) \leq f(R) n \|g\| L^{-1}
\]
by \eqref{73}, where~$f(R)$ is a suitable function of~$R$.  The
corresponding cycle~$M'$ has volume at most
\[
\left(n+\frac{n}{L} \right)\|g\| = n ||g|| (1+f(R)L^{-1}).
\]
Since~$M'$ admits a short map to~$X$, its volume is bounded below
by~$n$.  Thus,~$\tfrac{1}{n}M'$ is a cycle in~$X$ representing the
class~$g$, whose mass exceeds the mass of~$\tfrac{1}{n}M$ at most by
an arbitrarily small amount.  This yields a lower bound for~$\|g\|$
which is arbitrarily close to~$1$.  Note that similar arguments have
appeared in the work of I.~Babenko and his students \cite{Bab1, Bab2,
Ba3, BB, Ba06}, as well as the recent work of M.~Brunnbauer \cite{Bru,
Bru3}.
\end{proof}

\forget
\begin{question}
Is there a nonzero degree map from either~$\HP^4$ or~$\CP^8$, to
the~$16$-dimensional Cayley projective plane?
\end{question}
More specifically, we constructed a CW complex~$\CP^4 \cup e^3 \cup
e^7$ which is homotopy equivalent to~$\HP^2$.  Here the 6-dimensional
homology of the subcomplex~$\CP^4 \cup e^3$ turns out to be spherical,
which allows one to attach another cell.  The Cayley case is more
difficult.  Consider the 16-dimensional CW complex~$\HP^4 \cup e^5$
obtained from~$\HP^4$ by coning off a copy of~$\HP^1$.  Consider
the~$12$-dimensional homology of~$\HP^4 \cup e^5$.  Is its generator
spherical?  \forgotten

\section{The Cayley form and the Kraines form}
\label{eight}

The proof of the upper bound \eqref{15} for the optimal stable
$4$-systolic ratio depends on the calculation of the Wirtinger
constant~$\Wirt_4$ of~$\R^8$, \cf Corollary \ref{55}.

This section contains an explicit description \eqref{83c} of the
Cayley~$4$-form~$\cay$ in terms of a Euclidean basis.  The seven
self-dual forms appearing in the decomposition of~$\cay$ turn out to
have Lie-theoretic significance as a basis for a Cartan subalgebra of
the Lie algebra~$E_7$, discussed in detail in Section~\ref{nine}.  The
fact that~$\cay$ has unit comass constitutes the lower bound part of
the evaluation of the Wirtinger constant of~$\R^8$.  The upper bound
follows from the Lie-theoretic analysis of Section~\ref{nine}.

In more detail, let~$\{dx_1,dx_2,dx_3,dx_4\}$ denote the dual basis to
the standard real basis~$\{1,i,j,k\}$ for the quaternion
algebra~$\HH$.  Furthermore, let~$\{ dx_\ell, dx_{\ell'} \}$,
where~$\ell=1,\ldots,4$, be the dual basis for~$\HH^2$.  The three
symplectic forms~$\omega_i$,~$\omega_j$, and~$\omega_k$ on~$\HH^2$
defined by the three complex structures~$i,j,k$ are

\begin{equation}
\label{81b}
\begin{aligned}
\omega_i&=&dx_1\wedge dx_2+dx_3\wedge dx_4+dx_{1'}\wedge dx_{2'}+
dx_{3'}\wedge dx_{4'},\\ \omega_j&=&dx_1\wedge dx_3-dx_2\wedge
dx_4+dx_{1'}\wedge dx_{3'}- dx_{2'}\wedge dx_{4'},\\ \omega_k&=&
dx_1\wedge dx_4+dx_2\wedge dx_3+dx_{1'}\wedge dx_{4'}+ dx_{2'}\wedge
dx_{3'} .
\end{aligned}
\end{equation}
Let
\[
dx_{abcd}:=dx_a\wedge dx_b\wedge dx_c\wedge dx_d,
\]
where~$\{a,b,c,d\}\subset \{1,\ldots,4,1',\ldots,4'\}$. The
corresponding wedge squares satisfy
\begin{equation}
\begin{aligned}
\tfrac12\omega_i^2 &=
(dx_{1234}+dx_{1'2'3'4'})+(dx_{121'2'}+dx_{343'4'})+
(dx_{123'4'}+dx_{341'2'}), \\ \tfrac12\omega_j^2 &=
(dx_{1234}+dx_{1'2'3'4'})+(dx_{131'3'}+dx_{242'4'})-
(dx_{132'4'}+dx_{241'3'}), \\ \tfrac12\omega_k^2 &=
(dx_{1234}+dx_{1'2'3'4'})+(dx_{141'4'}+dx_{232'3'})+
(dx_{142'3'}+dx_{231'4'})
\label{cartanbasis}
\end{aligned}
\end{equation}
The seven distinct self-dual 4-forms appearing in
decomposition~\eqref{seven}~of the Kraines form, which are also
displayed in parentheses in~\eqref{cartanbasis}, form a basis of
a~$7$-dimensional abelian subalgebra~$\mathfrak h$ of the exceptional
real Lie algebra~$E_7$. In fact, the subalgebra that they generate is
a maximal abelian subalgebra of~$E_7$, as explained in
Section~\ref{nine}. The Cayley form
\[
\cay =\frac{1}{2}\left( \omega_i^2 + \omega_j^2 - \omega_k^2 \right)
\]
is the sum of the seven selfdual forms, with suitable signs, and
without multiplicities:
\begin{equation}
\label{83c}
\cay=e^{1234}+e^{1256}+e^{1278}+e^{1357}-e^{1467}-e^{1368}-e^{1458},
\end{equation}
where~$e^{abcd} = dx_{abcd}+ *dx_{abcd}$, while indices~$1'\ldots,4'$
are relabeled as~$5,\ldots,8$.
\begin{proposition}
\label{83b}
The Cayley form has unit comass.
\end{proposition}

\begin{proof}
R.~Harvey and H.~B.~Lawson \cite{HL} clarify the nature of the Cayley
form, as follows.  They realize the Cayley form as the real part of a
suitable multiple vector product on~$\R^8$ \cite[Lemma B.9(3),
p.~147]{HL}.  One can then calculate the comass of the Cayley form,
denoted~$\Phi$ in~\cite{HL}, as follows.  Let~$\zeta=x\wedge y\wedge z
\wedge w$ be a~$4$-tuple.  Then
\[
\Phi(\zeta)=\Re(x\times y\times z\times w)\leq | x\times y\times
z\times w | = |x\wedge y \wedge z \wedge w| ,
\]
and therefore~$\|\Phi\|=1$.  See also \cite{KSh} for an alternative
proof.
\end{proof}

By way of comparision, note that the square~$\eta= \tau^2$ of the
Kahler form~$\tau$ on~$\C^4$ satisfies~$\tfrac {|\eta|^2} {\|\eta\|^2}
=6$.  Meanwhile, the Cayley form yields a higher ratio, namely~$14$,
by Proposition \ref{83b}.  The Cayley form, denoted~$\omega_1$ in
\cite[p.~14]{DHM}, has unit comass, satisfies~$|\omega_1| ^2 =14$, and
is shown there to have the maximal ratio among all selfdual forms
on~$\R^8$.

The~$E_7$ viewpoint was not clarified in \cite{HL, DHM}.  Thus, the
``very nice seven-dimensional cross-section'' referred to in
\cite[p.~3, line~8]{DHM} and \cite[p.~12, line~5]{DHM}, is in fact a
Cartan subalgebra of~$E_7$, \cf Lemma~\ref{83}.

The calculation of~$\W_4$ results from combining Lemma~\ref{121} and
\cite{DHM}.  We will give a more transparent proof, using~$E_7$, in
the next section.

\section{$E_7$, Hunt's trick, and Wirtinger constant of~$\R^8$}
\label{nine}

To prove the upper bound of \eqref{15}, by
Proposition~\ref{upperbound}, we need to calculate the Wirtinger
constant of~$\R^8$.

\begin{proposition}
\label{45} We have~$\W_2=2$, while~$\W_4 = 14$.
\end{proposition}

\begin{proof}
By the Wirtinger inequality and Corollary~\ref{1214}, we obtain
the value~$\W_2= 2$.

To calculate the value of~$\W_4$, it remains to show that
no~$4$-form~$\omega$ on~$\R^8$ has a ratio~$|\omega|^2/\|\omega\|^2$
higher than~$14$.  By Lemma~\ref{121}, we can restrict attention to
selfdual forms.  We will decompose every such~$4$-form into the sum of
at most~$14$ simple (decomposable) forms with the aid of a particular
representation of a self-dual~$4$-form, stemming from an analysis of
the exceptional Lie algebra~$E_7$. Such a representation of a
self-dual~$4$-form was apparently first described explicitly by
L.~Antonyan \cite{An}, in the context of the study of~$\theta$-groups
by V.~Kac and E.~Vinberg \cite{GV} and E.~Vinberg and
A.~Elashvili~\cite{VE}.

We first recall the structure of the Lie algebra~$E_7$, following the
approach of J.~Adams~\cite{Ad}.  The Lie algebra~$E_7$ can be
decomposed as a direct sum
\begin{equation}
\label{82}
E_7=sl(8)\oplus\Lambda^4(8),
\end{equation}
\cf \cite[p.~76]{Ad}.  The Lie bracket on~$sl(8) \subset E_7$ is the
standard one.  The Lie bracket~$[a,x]$ of an element~$a\in sl(8)$ with
an element~$x\in \Lambda^4(8)$ is given by the standard action
of~$sl(8)$ on~$\Lambda^4(8)$.  Meanwhile, the Lie bracket of a pair of
elements~$x, y \in \Lambda^4(8)$ is defined as follows, \cf
\cite[p.~76, line~9]{Ad}:
\[
(a,[x,y])_{sl}=([a,x],y)_\Lambda^{\phantom{I}} .
\]
 The non-degenerate, but indefinite, inner product on~$sl(8)$
is given by
\[
(a,b)_{sl}=\trace ab.
\]
and the (non-degenerate, indefinite) inner product
on~$\Lambda^4(8)$ is given by
\[
(\alpha, \beta)_\Lambda^{\phantom{I}} \;d\!\vol=\alpha\wedge\beta,
\]
where~$d\!\vol$ is the volume form. If we complete this definition to
an inner product on~$E_7$ in which~$sl(8)$ and~$\Lambda^4(8)$ are
orthogonal, then the result is an invariant, non-degenerate,
indefinite inner product~$(\, ,\,)$ on~$E_7$ and the Killing form
is~$36(\, ,\,)$. See \cite [p. 78, ``Addendum"]{Ad}.

\begin{proposition}
\label{81}
In coordinates, the Lie bracket on~$\Lambda^4(8) \subset E_7$ can be
written as follows.  Let~$e_1,\ldots,e_8$ be a basis of
determinant~$1$.  Then
\begin{equation}
\label{lb}
\begin{aligned}
\left[ e_{r_1}e_{r_2}e_{r_3}e_{r_4}, e_{s_1}e_{s_2}e_{s_3}e_{s_4}
\right] &= 0 \quad \mbox{\rm if two or more~$r$'s equal~$s$'s}, \\
[e_1e_2e_3e_4, e_4e_5e_6e_7] &= e_4\otimes e_8^*, \\
[e_1e_2e_3e_4,e_5e_6e_7e_8] &={1\over 2}((e_1\otimes e_1^* +e_2\otimes e_2^*
+e_3\otimes e_3^* + e_4\otimes e_4^*)\\
&\quad-(e_5\otimes e_5^*+e_6\otimes e_6^*+e_7\otimes e_7^*+e_8\otimes e_8^*)).
\end{aligned}
\end{equation}
\end{proposition}

This is proved in \cite[p.~76]{Ad}.

The decomposition in \eqref{82} can be refined into the Cartan
decomposition of a Riemannian symmetric space for the group~${\bf
E_7}$, a non-compact form of ~${\bf E_7}/[SU(8)/\{\pm I\}]$, see
\cite[p.285]{Wo}. Recall that, in general, a Cartan decomposition of a
real Lie algebra consists of a maximal compact subalgebra, on which
the restriction of the Killing form is negative definite, and an
orthogonal positive definite complement. The Cartan decomposition
for the Riemannian symmetric space $SL(8,\R)/SO(8)$ is
\[
sl(8)=so(8)\oplus sym_0(8),
\]
where~$sym_0(8)$ is the set of~$8\times8$ traceless symmetric
matrices. The $SO(8)$ representation $\Lambda^4(8)$ is a direct
sum
\[
\Lambda^4(8)=\Lambda^4_+(8)\oplus \Lambda^4_-(8),
\]
where the subscripts~$+$ and~$-$ indicate ``selfdual" and
``anti-selfdual" forms, respectively.  Then the Cartan decomposition
for $E_7$ is given by
\begin{eqnarray*}
E_7&=&{\mathfrak k}\oplus {\mathfrak p}\\
{\mathfrak k}&=&so(8)\oplus\Lambda^4_-(8)\\
{\mathfrak p}&=&sym_0(8)\oplus \Lambda^4_+(8).
\end{eqnarray*}

One of the standard results in the theory of real reductive Lie groups
is the conjugacy of maximal abelian subalgebras of the noncompact
component~$\mathfrak p$ of the Cartan decomposition.  Here the term
``maximal abelian subalgebra'' refers to a subalgebra of~$\mathfrak p$
which is maximal with respect to the condition of being an abelian
subalgebra of~$E_7$, see \cite[\S 2.1.6, \S 2.3.4] {Wal}.  We will
apply the conjugacy condition inside a Lie subalgebra of $E_7$,
\[
E_7\supset \mathfrak g:= so(8) \oplus \Lambda^4_+(8)={\mathfrak
k}_0\oplus {\mathfrak p}_0
\]
and to a maximal abelian subalgebra $\mathfrak
h\subset\Lambda^4_+(8)={\mathfrak p}_0$ which contains the Cayley
form~$\cay$ \cite[Definition~10.5.1]{Jo0}.  The Cayley form is the
signed sum of~$7$ self-dual $4$-forms defining a basis of $\mathfrak
h$.  The exact expression for $\cay$ is given in \eqref{83c},
see~\cite{Br} and~\cite[equation~10.19]{Jo0}.

\begin{definition} Define the subspace~$\mathfrak h$ of
$\Lambda^4_+(8)$ as the span of the self-dual~$4$-forms
of \eqref{83c}, namely
\[
{\mathfrak h}= \R e^{1234} \oplus \R e^{1256} \oplus \R e^{1278}
\oplus \R e^{1357} \oplus \R e^{1467} \oplus \R e^{1368} \oplus \R
e^{1458}.
\]
\end{definition}
\begin{lemma}
\label{83} The subspace ~${\mathfrak h}$ is  a maximal abelian
subalgebra of~$\Lambda^4_+(8)$.
\end{lemma}

\begin{proof}
The bracket on $\mathfrak g$ is the restriction of the ~$E_7$ Lie
bracket  described in \cite[p.~76]{Ad} and Proposition~\ref{81}.
The bracket of two simple~$4$-forms vanishes whenever the forms
have a common~$dx_i\wedge dx_j$ factor, and it is easy to see that
this condition is satisfied for all the Lie brackets of pairs of
simple forms which occur in the Lie brackets of the seven self-dual
forms.  Since~$E_7$ is of rank~$7$, the dimension of a maximal abelian
subalgebra of~${\mathfrak p}$ is ~$7$,which gives and upper bound
on the dimension of an abelian subalgebra of $\Lambda^4_+(8)$.
\end{proof}

The following theorem shows that every self-dual~$4$-form
is conjugate by an element of~$SO(8)$ to an element
of~$\mathfrak h$, which completes the proof of
Proposition~\ref{45}.
\end{proof}

\begin{theorem}\cite[Theorem 8.6.1]{Wo}
\label{54}
Let $\mathfrak g=\mathfrak k \oplus \mathfrak p$ be the Cartan decomposition
associated to a Riemannian symmetric space $G/K$. Let $\mathfrak a$ and $\mathfrak a'$
be two maximal subalgebras of $\mathfrak p$. Then
\begin{enumerate}
\item there exist an element $X\in \mathfrak a$ whose centralizer in $\mathfrak p$
is just $\mathfrak a$,
\item  there is an element $k\in K$ such that $Ad(k){\mathfrak a}'=\mathfrak a$,
\item $\mathfrak p=\bigcup_{k\in K} Ad(k)\mathfrak a$.
\end{enumerate}
\end{theorem}

\begin{proof}[Partial proof of Theorem~\ref{54}]
 The proof of item (1) makes use of the compact dual symmetric space. 
  In the compact model the desired  element of the algebra 
is such  that the associated one parameter subgroup is dense in a maximal 
torus. For details of  the proof of (1) see \cite[page 253]{Wo}. We will 
prove (2) and (3), beginning with (3). The proof uses
an idea of G.~Hunt \cite{Hu}.

Let~$X\in {\mathfrak a}$ be the element whose existence is established in (1):
\[
{\mathfrak a}=\left\{ Y\in {\mathfrak p} \; \left| \; [Y,X]=0
\right. \right\}.
\]
Let $Z\in \mathfrak p$ be arbitrary. Consider the following function~$f$ on~$SO(8)$:
\[
f(k)=B(Ad(k)Z,X),
\]
where~$B(-,-)$ is the Killing form on~$\mathfrak g$.  Since~$SO(8)$ is
compact, the function attains a minimum at some point~$k$.  For
all~$W\in so(8)$, we have
\begin{eqnarray*}
0&=&\tfrac{d}{dt} |_{t=0}^{\phantom{I}} B(Ad(\exp(tW)k)Z,X)\\
&=&B([W,Ad(k)Z], X)\\&=& B(W, [Ad(k)Z,X])
\end{eqnarray*}
by the~$ad$-invariance of the Killing form.  Since the Killing form on
$so(8)$ is negative definite, it follows that~$[Ad(k)Z,X]=0$.
Thus~$Ad(k)Z\in {\mathfrak a}$, and $Z\in Ad(k^{-1})(\mathfrak a)$,
proving (3).

To prove (2) let $X'$ be an element whose centralizer in $\mathfrak p$ is
$\mathfrak a'$:
\[
{\mathfrak a}'=\{Y\in {\mathfrak p} \;|\; [Y,X']=0\}.
\]
We have just proved that there exists an element $k\in K$ such
that $[Ad(k)(X'), X]=0$; therefore, $Ad(k)(X')\in \mathfrak a$.
Thus $\mathfrak a$ centralizes $Ad(k)(X')$ and
$Ad(k^{-1}){\mathfrak a}$ centralizes $X'$; so
$Ad(k^{-1}){\mathfrak a}\subset {\mathfrak a}'$. Similarly,
$[Ad(k^{-1})(X), X']=0$ and $Ad(k){\mathfrak a}'\subset {\mathfrak
a}$. Thus $Ad(k){\mathfrak a}'= {\mathfrak a}$, concluding the
proof of (2).
\end{proof}
This completes the proof of Theorem~\ref{54} and hence of
Proposition ~\ref{45}.

\section{$b_4$-controlled surgery and systolic ratio}
\label{ten}

We will refer to an~$8$-manifold with exceptional~$\Spin(7)$ holonomy
as a {Joyce manifold}, \cf \cite{Jo0}.  Known examples of Joyce
manifolds have middle dimensional Betti number ranging from~$84$ into
the tens of thousands.  It is unknown whether or not a Joyce manifold
with~$b_4=1$ exists.  Yet no restrictions on~$b_4$ other than~$b_4\geq
1$ are known.  The obligatory cohomology class in question is
represented by a parallel Cayley~$4$-form~$\cayp$, \cf \eqref{83c},
representing a generator in the image of integer cohomology.

\begin{proposition}
\label{101}
A hypothetical Joyce manifold~$\hypo$ with unit middle Betti number
would necessarily have a systolic ratio of~$14$.
\end{proposition}

\begin{proof}
A generator of~$H^4(\hypo,\Z)_\R=\Z$ is represented by~$\cayp$. By
Poincar\'e duality, the square of the generator is the fundamental
cohomology class of~$\hypo$.  Thus, similarly to \eqref{P61d} and
\eqref{49}, we can write
\begin{equation}
\begin{aligned}
1 & = \int_{\hypo} \left| {\cayp}^{\wedge 2} \right|\; d\!\vol \\ & =
14 \left( \| \cayp \|_\infty \right)^2 \vol_{8}(\hypo),
\end{aligned}
\end{equation}
and the proposition follows by duality of comass and stable norm, as
in Gromov's calculation.
\end{proof}

The theorem below may give an idea of the difficulty involved in
evaluating the optimal ratio in the quaternionic case, as compared to
Pu's and Gromov's calculations.

\begin{theorem}
\label{14}
If there exists a Joyce manifold with~$b_4=1$, then the common value
of the middle dimensional optimal systolic ratio of~$\HP^2$ and
$\CP^4$ equals~$14$.  In particular, in neither case is the symmetric
metric optimal for the systolic ratio.
\end{theorem}

We introduce a convenient term in the context of surgery on an
$8$-dimensional manifold~$M$.

\begin{definition}
A~$b_4$-{\em controlled surgery\/} is a surgery which induces an
isomorphism of the~$4$-Jacobi torus \eqref{13}.
\end{definition}

In particular, such a surgery does not alter the middle dimensional
Betti number~$b_4(M)$.  It was shown in Section~\ref{sevenb} that such
a surgery does not alter the stable~$4$-systolic ratio.

\begin{proposition}
\label{22sur}
Every simply connected spin~$8$-manifold~$M$ satisfying~$b_4(M)=1$
admits a sequence of~$b_4$-controlled surgeries, resulting in
a~$2$-connected manifold, denoted~$\rhp2$, with the rational
cohomology ring of the quaternionic projective
plane:~$H^*(\rhp2,\Q)=H^*(\HP^2,\Q)$.
\end{proposition}

\begin{proof}
We choose a system of generators~$(g_i)$ for~$H_2(M,\Z)$.  By the
Hurewicz theorem, each~$g_i$ can be represented by an imbedded
$2$-sphere~$S_i \subset M$.  The spin condition implies the triviality
of the normal bundle of each~$S_i$.  We can therefore perform
successive surgeries along each~$S_i$ to remove~$2$-dimensional
homology, resulting in a~$2$-connected manifold~$M'$.  Clearly,
$b_4(M')=1$, while the third Betti number may have changed during the
surgeries.

Similarly, we choose a system of~$3$-spheres representing a basis for
$H_3(M',\Q)$.  The normal bundles are automatically trivial, and
surgeries along the~$3$-spheres reduce the~$b_3$ to zero without
altering~$b_4$, resulting in a manifold~$\rhp2$ with the rational
cohomology of the quaternionic projective plane by Poincar\'e duality.
\end{proof}

\begin{corollary}
\label{104}
A Joyce manifold with~$b_4=1$ admits a sequence of~$b_4$-controlled
surgeries which produce a manifold~$\rhp2$ which has the rational
cohomology of~$\HP^2$.
\end{corollary}

\begin{proof}
Manifolds with~$\Spin(7)$ holonomy are simply connected and
spin by \cite[Theorem 10.6.8, p.~261]{Jo0}, and we apply
Proposition~\ref{22sur}.
\end{proof}

Note that the ``cylinder'' of a surgery transforming~$X$ to~$Y$ is
homotopy equivalent to a complex~$W$ obtained from~$X$ by attaching a
cell.  Thus the inclusion of~$Y$ as the other end of the cylinder
defines a map~$Y\to W$ indicing an isomorphism of the Jacobi torus
$J_4$.  Applying the pullback techniques of Section~\ref{sevenb}, we
conclude that~$\SR_4(X)=\SR_4(Y)$.  An interesting related
axiomatisation (in the case of 1-systoles) is proposed in \cite{Bru}.

\begin{proposition}
\label{Sh}
A manifold~$\rhp2$ with the rational cohomology of the quaternionic
projective plane admits a nonzero degree map~$\HP^2 \to \rhp2$ from
$\HP^2$.
\end{proposition}

\begin{proof}
The fact that~$\HP^2$ has a map of nonzero degree to a manifold with
its rational cohomology algebra, follows from the formality of the
space combined with the theorem of H.~Shiga \cite{Sh}.  Namely, the
theorem gives enough self maps of any formal space to build its
rational homotopy type by iterated mapping cylinders.  Hence~$\HP^2$
admits a map to the rationalisation of~$\rhp2$.  By compactness, the
image of the map lies in a finite piece of the iterated space.  The
finite piece admits a retraction to~$\rhp2$ itself.  This produces the
desired map.
\end{proof}

\begin{corollary}
\label{107}
A manifold~$\rhp2$ with the rational cohomology of the quaternionic
projective plane satisfies~$\SR_4(\HP^2) \geq \SR_4(\rhp2)$.
\end{corollary}

\begin{proof}
Let~$d^2$ be the degree of the map.  We then constuct suitable metrics
on the quaternionic projective plane by pullback.  The argument is
similar to that of Section~\ref{sevenb} and relies upon the existence
of~$d$-motonone maps, \ie maps such that the preimage of a
path-connected set has at most~$d$ path connected components, see
\cite{Wr, Bru, Bru3}.  In more detail, we have $\vol(\HP^2) = d^2
\vol(\rhp2)$.  Meanwhile, the induced homomorphism in $H_4$ is
multiplication by $d$.  Since the stable norm is by definition
multiplicative.  Hence $\stsys_2(\HP^2) \geq d \stsys_s(\rhp2)$,
proving the corollary.
\end{proof}

\begin{remark}
\label{107b}
A referee asked whether the map in Proposition~\ref{Sh} can be taken
to be of degree~$1$.  Whereas in general this is not the case, it
turns out that in the absence of torsion, degree~$576$ is sufficient,
as shown in Section~\ref{eleven}.
\end{remark}

\begin{proof}[Proof of Theorem \ref{14}]
A Joyce manifold has systolic ratio of~$14$ by Proposition~\ref{101}.
By Corollary~\ref{104}, the manifold~$\rhp2$ must also
satisfy~$\SR_4(\rhp2)=14$.  Finally, Corollary~\ref{107} implies
that~$\SR_4(\HP^2)=14$, as well.
\end{proof}

\section{Hopf invariant, Whitehead product, and systolic ratio}
\label{eleven}

This section answers a question referred to in Remark~\ref{107b}.
S.~Smale as well as J.~Eells and N.~Kuiper \cite{EK} proved that every
manifold which is a \fake~$\HP^2$, is homotopy equivalent to~$S^4
\cup_h e^8$, where the attaching map lies in a class
\begin{equation}
\label{111b}
[h]\in \pi_7(S^4)=\Z + \Z_{12}
\end{equation}
which is an infinite generator.

Let~$m\geq 2$ be an even integer.  Let~$e\in \pi_m(S^m)$ be the
fundamental class.  Let~$q\geq 1$, and consider a self map of~$S^m$ of
degree~$q$.  Let
\begin{equation}
\label{71b}
\phi_q: \pi_{2m-1}(S^m) \to \pi_{2m-1}(S^m)
\end{equation}
be the induced homomorphism.  The following result is immediate from
standard properties of Whitehead products~$[\;,\;]$.

\begin{lemma}
\label{111}
The class~$[e,e]\in \pi_{2m-1}(S^m)$ satisfies~$\phi_q([e,e])=q^2
[e,e]$.
\end{lemma}

Given an element~$x\in \pi_{2m-1}(S^m)$, we can write
\begin{equation}
\label{71z}
2x=s+H(x)[e,e],
\end{equation}
where~$s$ is torsion, and~$H(x)$ is the Hopf invariant of~$x$.  Note
that if~$x$ is the class represented by the Hopf fibration, then~$s$
is a generator of the torsion subgroup.  In particular, the
class~$[e,e]$ is primitive (\ie not twice another class) in the
quaternionic case, unlike the complex case.

We have the following formula for the map \eqref{71b}, \cf B.~Eckmann
\cite{Ec} and G.~Whitehead \cite[p.~537]{Wh}:
\begin{equation}
\label{72}
\phi_q(x)=qx+ {q \choose 2} H(x) [e,e].
\end{equation}

\begin{lemma}
\label{73z}
For all~$x\in \pi_{7}(S^4)$, if~$q$ is a multiple of~$24$, then
\[
\phi_q(x)=q^2x = \tfrac{q^2}{2}[e,e].
\]
\end{lemma}

\begin{proof}
Let~$a$ be the attaching map of the true~$\HP^2$.  By \eqref{71z} and
\eqref{111b}, the multiple~$qa$ (and hence~$q^2a$) is proportional
to~$[e,e]$.  Therefore by~\eqref{72}, the image~$\phi_q(a)$ is also
proportional to~$[e,e]$.  Thus,~$\phi_q(a)$ is proportional to every
infinite generator~$x$ by Lemma~\ref{111}, proving the lemma.
\end{proof}

\begin{theorem}
\label{576}
Any \fake~$\HP^2$ admits a continuous map of degree~$576$ from the
true~$\HP^2$.
\end{theorem}

\begin{proof}
By Lemma~\ref{73z}, a self-map of~$S^4$ of degree a multiple of~$24$,
necessarily sends the attaching map of the true~$\HP^2$, to a class
proportional to the attaching map of the \fake{} one.  Hence the map
can be extended over the entire 8-manifold.
\end{proof}

\section*{Acknowledgements}
We are grateful to D.~Alekseevsky, I.~Babenko, R.~Bryant,
A.~Elashvili, D.~Joyce, V.~Kac, C.~LeBrun, and F.~Morgan for helpful
discussions.

\vfill\eject

\end{document}